\documentclass[11pt]{amsart}
\usepackage{amssymb,enumerate,enumitem,commath,xcolor, graphicx,mathtools,bbm,dsfont,hyperref}
\usepackage[sorting=none]{biblatex}
\calclayout

\newtheorem{thm}{Theorem}[section]

\newtheorem{cor}[thm]{Corollary}

\theoremstyle{definition}
\newtheorem{dfn}[thm]{Definition}

\theoremstyle{remark}

\numberwithin{equation}{section}

\newcommand{\R}{\mathbb{R}}  

\newcommand{\N}{\mathbb{N}}

\begin{document}

\title{Yet Another Quantitative Harris Theorem }
\author{Christopher DuPre}
\address{Department of Mathematics, Georgia Institute of Technology, 
Atlanta, GA 30332}
\email{cdupre3@gatech.edu}


\begin{abstract}
In this paper we develop a quantitative Harris theorem with effective control over the constants. A benefit of our methodology is the decoupling of the small set and Lyapunov-Foster Drift conditions. Our methodology allows any small set and any set in the Lyapunov-Foster condition as long as the second satisfies a so-called ``quantitative petiteness" condition. The theorem relies on a novel proof of a quantitative Kendall-type theorem
which is inspired by the techniques of Markov Chains on general state spaces. We give an application of the technique to the Markov chain approximation of mixing processes.
\end{abstract}
\maketitle
\section{Introduction}
In this paper, we will present a novel proof of a quantitative mixing theorem for Markov chains on general state spaces. Such limiting statements for the evolution of Markov chains on general state spaces are collectively referred to as Harris Theorems in honor of the pioneering work of T.E. Harris in the fifties \cite{harris1960existence}. We follow the work of Meyn and Tweedie in \cite{meyn2012markov} by extending the ideas of Doeblin in the creation of a ``small" set to a general state space. This is then combined with a technique known as Numellin splitting \cite{nummelin1978splitting} to reduce the problem to a statement of renewal theory. Our particular approach relies on mixing with respect to a weigthed norm. This technique was introduced in \cite{hordijk1992ergodicity}.
\par In the mixing case, this approach relies on a key theorem pioneered by Kendall \cite{kendall1959unitary}. Most proofs in the literature seem to rely on techniques from complex analysis \cite{meyn2012markov,baxendale2005renewal}. While these proofs are beautiful in their own right, we present a new proof of a major implication in Kendall's theorem which we need for our estimates. The proof relies on techniques generated for general state space Markov chains and a coupling approach. As far as we can tell, this seems to be the first such coupling proof in the literature.
\par This technique has been successful in the creation of both qualitative \cite{meyn2012markov,rosenthal1995minorization} and quantitative\cite{baxendale2005renewal, WitoldBednorz2013,jiang2020coupling} estimates on mixing rates. There have also been several impressive applications of these ideas to rigorous or ``honest" error bounds for Markov Chain Monte Carlo algorithms \cite{latuszynski2011rigorous,jones2001honest,tierney1994markov}. 
\par The motivation for this work came from considering uniform mixing rates a parametric family of random diffeomorphisms on a compact state space as an extension of \cite{blumenthal2022exponential}. In our case, the alignment of a small set and the sub-level set was unnatural and so we created these tools. We hope they may be of independent interest.
\section{Setting and Main Theorem}
In this paper, we will be considering discrete-time Markov Chains on general state spaces. By this we mean that we have a measurable space $\left(X,\mathcal{B}\right)$ and a kernel $P(x,\cdot)$ defined such that 
\begin{enumerate}
    \item For every $x\in X$, $P(x,\cdot)$ is a probability measure
    \item For every $A\in \mathcal{B}$, $ P(x,A)$ is a non-negative measurable function.
\end{enumerate}
This kernel encodes the one-step transition probabilities of our Markov process. The n-step transition probabilities are then defined inductively via the Chapman-Kolmogorov equations.
$$P^n(x,\cdot) := \int_X P^{n-1}(x,dy) P(y,\cdot). $$
We can define an action of our kernel on measurable functions $f$ and measures $\lambda$ as follows.
\begin{align*}
   P^n(\lambda,\cdot) := \lambda P^n (\cdot) := \int_X \lambda(dy) P^N(y,\cdot) && P^n(x,f) := P^n f(x) := \int_X P^N(x,dy) f(y).
\end{align*}
Given any initial measure $\lambda$, we then define the Markov chain $\{\Phi^n\}_{n\in\N}$ as the stochastic process such that for every increasing sequence of indices $i_1 < i_2 < ... i_n$ we have the following condition on the law.
\begin{gather*}
    \mathbb{P}\left(\Phi^{i_1} \in A_1, ... \Phi^{i_n} \in A_n \right) = \int_{A_1}\int_{A_2} ... \int_{A_{n-1}} P^{i_1}(\lambda,dy) P^{i_2-i_1}(y,dz)... P^{i_n-i_{n-1}}(w,A_n).
\end{gather*}
We will use the notation $\mathbb{P}_{x}$ to refer the law induced by the Markov chain with initial measure $\delta_x$. We will use $\mathbb{E}_{x}$ to denote integrating with respect to this law. For any $C\in \mathcal{B}$, let $\tau_C$ denote the first hitting time.
\par Our goal will be to obtain quantitative estimates on the mixing rate for our Markov chain. If $f:X\to \R$ is a measurable function on our space, we will use the notation of $f(\Phi^n)$ to mean the composition of the n-th step in our Markov chain and the function $f.$ An important construction of our work will be the so called $V$-weighted norm. For every measurable $V:X\to [1,\infty)$, we define a norm on measurable functions from $X$ as follows.
$$\|\phi\|_V := \sup_{x\in X} \frac{\abs{\phi(x)}}{V(x)}. $$
In particular, every bounded measurable function is also $\|\cdot\|_{V}$ bounded.
\begin{thm}[Quantitative Harris Theorem, \cite{meyn2012markov}]
\label{MainTheorem}
Let $P$ be a Markov transition kernel defined on a measurable space $(X,\mathcal{B})$ and assume the following:
\begin{enumerate}[label=(\alph*)]
    \item (Small set) There exists a measurable set $U$, a $\delta > 0$, and a probability measure $\mu$ on $X$ such that $\mu(U) = 1$ and 
    $$P(x,B) \geq \delta \mu(B) \text{ for all }x\in U,B \in \mathcal{B}. $$
    \item (Drift condition) There exists a measurable function $V:X\to [1,\infty)$, a constant $0<\lambda<1$, a constant $b<\infty$, and a set $C\in \mathcal{B}$ such that 
    $$PV \leq \lambda V + b\mathbbm{1}_{C}.  $$
    \item (Quantitative Petiteness) There exists an $N_0 \in \N$ and a $c>0$ such that 
    $$\underset{x \in C }{\inf} \mathbb{P}_{x}\left\{ \tau_U \leq N_0 \right\} \geq c. $$
\end{enumerate}
Let $M_U = \underset{x\in U}{\sup}\ V(x)$ and let $M_C = \underset{x\in C}{\sup}\ V(x).$
Then there exists an invariant measure $\pi$ and a pair of effective constants $D<\infty,0<\gamma<1$ which depend only (and explicitly) on $\delta,\lambda,b,N_0,c, M_U,M_C$ such that 
$$\abs{P^n\phi(x) - \int \phi d\pi} \leq DV(x)\|\phi\|_{V} \gamma^n. $$
\end{thm}
Note that several qualifications for this theorem are unnecessary. $U$ can be a small set for the n-step chain for any n. The measure $\mu$ need not give $U$ a positive mass, but only needs to satisfy an aperiodicity requirement. In general, the proof can be expanded to the reader's exhaustion. We take the above requirements for the sake of clarity.
\par A key technique in the study of ergodicity in general state Markov chains is the reduction to a problem in renewal theory \cite{meyn2012markov}. As such, we have developed a new proof of a key result related to Kendall's theorem. While Kendall's theorem is typically richer in structure than our result, this result will be all we need for our purposes. 
\begin{thm}
\label{MainKendallTheorem}
Let $p(n)$ be an increment distribution such that $p(1) \geq \beta > 0$ and let $u(n)$ be the associated renewal sequence with 0 delay. If there exists a $r>1$ such that
$$\sum_{n=1}^\infty p(n) r^n  \leq B $$
then there exists an invariant distribution for the forward recurrence time chain $e(n)$, a constant $\pi(1) = e*u(n)$, and an effective bound $\rho$ depending only (and explicitly) on $r,\beta,B$ such that for any $r_2 \leq \rho$ we can find an effective bound $L$ depending only (and explicitly) on $r_2,B,\beta, \rho$ such that
$$\sum_{n=1}^\infty \abs{u(n)-\pi(1)}r_2^n \leq L. $$
\end{thm}
This is again restricted by clarity instead of necessity. In fact the condition $p(n)>0$ and $p(m) >0$ for $n,m$ co-prime is sufficient. This is clear from our proof and does not essentially affect finding the necessary quantitative bounds. Our proof appears to be novel in that it appears to be the first coupling proof of this implication of Kendall's theorem in the literature.
\section*{Appendix A. Main Tools}
In this section we introduce the main tools which we will use in the statement of both theorems. The first tool gives explicit bounds on the moment generating function of $r^{\tau_C}$ where $C$ is the set in our Lyapunov Foster drift condition.
The result is effectively a consequence of Dynkin's Formula.
\begin{thm}[Theorem 15.2.5 in \cite{meyn2012markov}]
\label{CReturnTimes}
If $V:X\to [1,\infty)$ satisfies the Lyapunov-Foster Drift condition
$$PV \leq \lambda V + b \mathbbm{1}_{C} $$
for $\lambda < 1, b < \infty, C\in \mathcal{B}$, then for any $r\in \left(1,\frac{1}{\lambda}\right)$
$$\frac{\mathbb{E}_{x}[r^{\tau_C}]-1}{r-1} \leq \mathbb{E}_{x}\left[\sum_{k=0}^{\tau_{C}-1} V(\Phi^k)r^k\right] \leq \frac{1+rb}{1-\lambda r}V(x).$$
\end{thm}
This allows us to use bounds on the function in our Lyapunov-Foster drift condition to obtain explicit bounds on the first hitting times to the set in our drift condition. However, we will need to transfer such bounds from set to set. For this, we will use the following result which follows from a close reading of Theorem 15.2.1 in \cite{meyn2012markov}.
\begin{cor}
\label{MainCor}
Suppose there exists an $r_1 > 1$ such that 
\begin{equation*}
\label{First Bound}
    \sup_{x\in C} \mathbb{E}_{x}\left[ \sum_{k=0}^{\tau_C-1} V(\Phi^{k})r_1^k \right] \leq M_0 <\infty. 
\end{equation*}
Further, assume that for $B\in \mathcal{B}$ there exists an $N_0\geq 1$ and a $c>0$ such  that
\begin{equation*}
\label{ Second bound}
    \inf_{x\in C} \mathbb{P}_{x}\left\{ \tau_B \leq N_0 \right\}\geq c. 
\end{equation*}
Then, there exists an effective $\rho$ depending only (and explicitly) on $r_1, M_0,c, N_0$ such that for all $r_2 < \rho$ we can find an effective constant $D$ depending only (and explicitly) on $r_1, M_0,c, N_0, r_2$ such that
\begin{gather*}
    \frac{\mathbb{E}_{x}\left[r_2^{\tau_B} \right]-1}{r_2-1} \leq \mathbb{E}_{x}\left[\sum_{k=0}^{\tau_B -1} V(\Phi^{k})r_2^k \right] \leq D V(x)
\end{gather*}
\end{cor}
\section*{Appendix B. Proof of Theorem \ref{MainKendallTheorem}}
We begin by introducing the language of renewal theory. 
Let $p(n): \N \to \R_{\geq 0}$ and \newline $d(n): \N\cup\{0\} \to \R_{\geq 0}$ be sequences of non-negative numbers such that $$\sum_{n=1}^\infty p(n) = \sum_{n=0}^\infty d(n) = 1.$$
We will refer to $p(n)$ as the increment distribution and $d(n)$ as the delay distribution.
\par We let $\{Y_i\}_{i=1}^\infty$ be independent copies of discrete random variables with the law $\mathbb{P}\left\{ Y_i = n \right\} = p(n)$ and let $Y_0$ denote a random variable with values in $\N\cup \{0\}$ with law $d(n)$ which is independent of all other $Y_i$. We then define $S_n^d = \sum_{i=0}^n Y_i$ and 
$$V_n^d = \begin{cases} 1 & S_m^d = n \text{ for some m}\\ 0 & \text{else} \end{cases}. $$
The associated renewal process is
$$ v_d(n) := \mathbb{P}\left\{V_n^d = 1 \right\} = \sum_{m=0}^\infty d*p^{*m}(n). $$
Of special importance are the renewal sequences with $d(0) = 1$. We denote such a sequence with the symbol $u(n) := v_0(n).$ We can the write a renewal sequence with any delay in terms of the zero delay sequence as $v_d(n) = d*u(n).$ For more details into renewal sequences we refer the reader to \cite{lindvall2002lectures}. 
\par Notice that $\sum_{n=1}^\infty p(n) r^n < \infty$ implies that $\sum_{n=1}^\infty p(n)n < \infty $ as in particular we have that $\frac{r}{r-1}r^n \geq n$ for all $r>1$. Let $e(n)$ denote a delay distribution with the law
$$e(n) :=  \frac{\sum_{j={n+1}}^\infty  p(j)}{\sum_{n=1}^\infty p(n)n} = \frac{1- \sum_{j=1}^{n} p(j)  }{\sum_{n=1}^\infty p(n)n}, n \in \N\cup\{0\}. $$
This can be easily seen to be the invariant probability measure of the forward recurrence time chain (see \cite{meyn2012markov,lindvall2002lectures} for more details). 
In particular, we have that
$$\pi(1) := \mathbb{P}\left\{V_n^e = 1 \right\}= e*u(n) = \frac{1}{\sum_{n=1}^\infty np(n)}$$
By definition of the renewal sequences, we have the reduction that 
$$\abs{u(n) - \pi(1)} = \abs{\mathbb{P}\left\{V^0(n) = 1 \right\}- \mathbb{P}\left\{ V^e(n) = 1\right\}}.$$ 
To control this difference, we introduce a coupling. 
\begin{dfn}
The \textbf{coupling time} of $V^0(n), V^e(n)$ is defined as the stopping time
$$T_{0,e} = \min\left\{n \in \N\cup\{0\} \bigg \vert V^0(n) = V^e(n) = 1 \right\}.$$
\end{dfn}
The coupling time represents the first simultaneous renewal of the two sequences. Note that after a renewal, the distribution of the next renewal time is uniquely determined by the increment sequence. Thus we will see that the coupling time allows us to control the difference between the two sequences after any given time.
\par To formalize this, let us define the random variables 
$$V^{0,e}(n) = \begin{cases} 
V^0(n) & n < T_{0,e}\\
V^e(n) & n \geq T_{0,e}.
\end{cases}$$
Notice that for any $n$ (due to the distributions being identical after the first simultaneous renewal) we have that $$\mathbb{P}\left(V^{0,e}(n) = 1\right) = \mathbb{P}\left(V^{0}(n) = 1 \right).$$
We then have that
\begin{align*}
    \abs{u(n) - \pi(1)} &= \abs{\mathbb{P}\left\{V^0(n) = 1 \right\}- \mathbb{P}\left\{ V^e(n) = 1\right\}} = \abs{\mathbb{P}\left\{V^{0,e}(n) = 1 \right\}- \mathbb{P}\left\{ V^e(n) = 1\right\}}\\
    &= \bigg \vert \mathbb{P}\left\{V^{0}(n) = 1, T_{0,e} > n \right\} + \mathbb{P}\left\{V^{e}(n) = 1, T_{0,e} \leq  n \right\}\\
    &\ \ \ \  - \mathbb{P}\left\{ V^e(n) = 1, T_{0,e} > n\right\} - \mathbb{P}\left\{ V^e(n) = 1, T_{0,e} \leq n\right\}\bigg \vert\\
    & \leq \mathbb{P}\{T_{0,e} > n \}.
\end{align*}
Thus, multiplying by the appropriate power of $r$ to be chosen later and summing, we have that
\begin{equation}
\label{Kendally}
    \sum_{n=1}^\infty \abs{u(n) - \pi(1)}r^n \leq \sum_{n=1}^\infty r^n \mathbb{P}\{T_{0,e} > n\} \leq \frac{r\mathbb{E}[r^{T_{0,e}}]}{r-1} = \frac{r}{r-1}\sum_{n=1}^\infty e(n)\mathbb{E}[r^{T_{0,n}}].
\end{equation}
Thus it suffices for us to obtain geometrically bounded tails for the first coupling time of a zero delay distribution and a constant delay distribution.
\par Inspired by methods in the theory of Markov chains on general state spaces (such as Theorem \ref{CReturnTimes} and Corollary \ref{MainCor}), we note that the boundedness of the moment generating function of a return time is related to the existence of a Lypaunov function. Thus we would like to cast our problem as a Markov chain on some state space, and then find an appropriate Lyapunov function.
\par We first consider the so-called forward recurrence time chain. This is a Markov chain defined on the state space $\N$ with the discrete $\sigma$-algebra with the transition kernel 
$$P(n,n-1) = 1 \text{ if }n > 1, P(1,m) = p(m). $$
Considering the evolution of our two renewal sequences independently, we can then imagine the state of our system a a point in $\N\times\N$. Initial delays of $0,n$ corresponds to an initial point of $(1,n+1)$. The system then induces a Markov chain on $\N\times\N$ with kernel defined as follows
\begin{align*}
    P\left((a,b),(a-1,b-1)\right) = 1 & \text{ if }a,b\neq 1\\
    P\left((1,b),(n,b-1)\right) = p(n) & \text{ if }a=1,b\neq 1\\
    P\left((a,1),(a-1,n)\right) = p(n) & \text{ if }a\neq 1,b=1\\
    P\left((1,1),(n,m)\right) = p(n)p(m) & \text{ if }a,b=1.
\end{align*}
A simultaneous renewal is equivalent to this Markov chain reaching $(1,1)$. Thus we have that $\mathbb{E}[r^{T_{0,n}}] = \mathbb{E}_{1,n+1}\left[r^{\tau_{1,1}}\right].$
\par Now consider the function $V(a,b) = \frac{r^{a-1}+r^{b-1}}{2}$. By direct computation, we have the following relations.
\begin{align*}
    P\left((a,b), V\right) & = \frac{1}{r}V(a,b) & \text{if }a,b\neq 1\\
    P\left((1,b),V\right) & = \frac{1}{2r}\left(\sum_{n=1}^\infty p(n) r^n + r^{b-1}\right) & a=1,b\neq 1\\
    P\left((a,1),V\right) & = \frac{1}{2r}\left( r^{a-1} + \sum_{n=1}^\infty p(n) r^n \right)& a=1,b\neq 1\\
    P\left((1,1),V\right) & = \frac{1}{r}\sum_{n=1}^\infty p(n) r^n & a=1,b\neq 1.
\end{align*}
We have that $ \sum_{n=1}^\infty p(n) r^n < \infty$ so at the very least each of these is well-defined. We wish the set $C$ that we will construct in the Lyapunov function to be finite for reasons we shall get to about controlling access from our set $C$ to $(1,1)$. Thus, the point $(1,1)$ is of no concern but each of the rays $(1,n),(n,1)$ are seemingly problematic. However, we have that
$$\frac{P\left((a,1),V\right)}{V(a,1)} =
\frac{r^{a-1}+ \sum_{n=1}^\infty p(n) r^n}{r^a + r}. $$
Using $\sum_{n=1}^\infty p(n)r^n \geq r$ for all $r>1$, we can rearrange the condition that this fraction is less than $\frac{1}{r} < \eta < 1$ as follows.
\begin{align*}
    \frac{P\left((a,1),V\right)}{V(a,1)} \leq \eta & \iff \frac{r^{a-1}+ \sum_{n=1}^\infty p(n) r^n}{r^a + r} \leq \eta\\
    &\iff \sum_{n=1}^\infty p(n) r^n -\eta r \leq r^{a-1}\left(\eta r-1\right)\\
    &\iff a \geq \frac{1}{\ln(r)}\ln\left(\frac{\sum_{n=1}^\infty p(n) r^n-\eta r}{\eta r-1}\right)+1.
\end{align*}
The same calculation applies in the situation when $a=1,b\neq 0$. Thus we have the Lyapunov-Foster Drift condition as follows.
\begin{align*}
    M &:= \left\lceil\frac{1}{\ln(r)}\ln\left(\frac{\sum_{n=1}^\infty p(n) r^n-\eta r}{\eta r-1}\right)+1 \right\rceil\\
    C &:= \{(1,b) \in \N^2 \text{ such that }b\leq M \}\cup \{(a,1) \in \N^2 \text{ such that }a\leq M \}\\
    PV&\leq \eta V + \frac{\max\left\{2\sum_{n=1}^\infty p(n) r^n,\sum_{n=1}^\infty p(n) r^n + r^{M-1} \right\}}{2r}\mathbbm{1}_{C}.
\end{align*}
In particular, using Theorem \eqref{CReturnTimes} we have the following bound. $$\sup_{x\in C} \mathbb{E}_{x}\left[\sum_{k=0}^{\tau_{C}-1} V(\Phi^k)r^k\right] \leq \frac{2+\max\left\{2\sum_{n=1}^\infty p(n) r^n,\sum_{n=1}^\infty p(n) r^n + r^{M-1} \right\}}{4(1-\eta r)}(r^{M-1} + 1) .$$
We now wish to transfer this bound from a bound on $\tau_C$ to a bound on $\tau_{1,1}$. From the construction of our Markov chain, it is clear that 
$ \underset{x\in C}{\inf}\mathbb{P}_{x}\{\tau_{1,1} \leq M \} \geq p(1)^M. $
This bound will act as a ``quantitative petiteness" condition as in the statement of our Harris theorem which allows us to transfer bounds on $\tau_C$ to $\tau_{1,1}$ due to Corollary \ref{MainCor}. We thus now have an effective constant $D$  and an effective new $r$ with a bound 
$$\mathbb{E}[r^{T_{0,n}}] = \mathbb{E}_{1,n+1}[r^{\tau_{1,1}}] \leq DV(1,n+1) = \frac{D(1+r^n)}{2} \leq \frac{D}{2}\left(\frac{1}{r-1} + 1 \right)r^n = \frac{Dr}{2(r-1)}r^n.$$ 
We have used the fact that $\frac{r^n}{r-1} > 1$ for all $n\geq 1$ for all $r> 1$.
Summing these up then gives the following.
\begin{align*}
    \mathbb{E}[r^{T_{0,e}}] &= \sum_{n=0}^\infty e(n)\mathbb{E}[r^{T_{0,n}}] \leq \frac{Dr}{2(r-1)}\sum_{n=0}^\infty e(n)r^n = \frac{Dr}{2(r-1)}\frac{1}{\sum_{n=1}^\infty np(n)} \left[\sum_{n=0}^\infty \sum_{j=n+1}^\infty p(j) r^n\right]\\
    & \leq \frac{Dr}{2(r-1)}\left[\sum_{j=1}^\infty \sum_{n=0}^{j-1} p(j)r^n\right] \leq \frac{Dr}{2(r-1)^2}\left[\sum_{j=1}^\infty p(j)r^j\right].
\end{align*}
Plugging this back into \eqref{Kendally} gives the result. 
\section*{Appendix C. Proof of Theorem \ref{MainTheorem}}
First notice that by scaling it is sufficient to consider $\phi$ such that $\|\phi\|_{V} \leq 1.$ We adopt a coupling based approach. We wish to control two separate quantities
\begin{enumerate}
    \item Geometric tails on $\tau_U$
    \item The construction of a so called ``atom" which can translate our problem into a renewal theory problem.
\end{enumerate}
Our main tools in Theorem \ref{CReturnTimes} and Corollary \ref{MainCor} give an effective way to achieve the first result.
Combining $M_U, M_C , c, N_0$ in the theorem, we get geometric tails for the first hitting time for $U$ from an arbitrary initial point.
\par We now return to the second of our goals. We recall the definition of an atom. 
\begin{dfn}
Let $P(x,dy)$ be a transition kernel on a measurable space $(X,\mathcal{B})$. A set $\alpha \in \mathcal{B}$ is called an \textbf{atom} if there exists a probability measure $\nu$ such that that
$$P(x,B) = \nu(B) \text{ for all } x\in \alpha. $$
\end{dfn}
We will reserve the symbol $\alpha$ for an atom. When $\alpha$ is an atom, we will write $P(\alpha,B)$ to mean $P(x,B)$ for any $x\in \alpha$ or equivalently $\nu(B)$. If a stopping time refers to an atom, we are referring to a Markov chain with an atom. We will construct such a chain in the next paragraph using the technique of Nummelin Splitting \cite{nummelin1978splitting}. We will then transfer the geometric tails to the newly created atom. 
\par We then split the chain at $U$ with a $\delta/2$ splitting. What we mean by this is that we consider the extended space $X\times\{0,1\}$ equipped with the product $\sigma$-algebra. We will use the notation that $ A_0 = A\times \{0\}, A_1 = A \times \{1\}$ $A\in \mathcal{B}$, $x_{i} = x\times\{i\}$ where $x\in X$. Given any measure $\lambda$ on $(X,\mathcal{B})$, we will construct a measure on this new space as 
\begin{align*}
    \lambda^*(A_0) &=  \left(1-\frac{\delta}{2}\right)\lambda(A\cap U) + \lambda(A\cap U^c) \\
    \lambda^*(A_1) &= \frac{\delta}{2}\lambda(A\cap U).
\end{align*}
Likewise, given any function $V:X\to \R$, we can define the function $\hat{V}:X\times\{0,1\} \to \R$ on the split chain as $\hat{V}(x_i) := V(x).$ We now are able to extend measures into this new space, so now we turn to extending transition probabilities. In this way, we will be able to construct a Markov Chain on this new space from our original chain. 
\par We define the split kernel $\hat{P}(x_i, A)$ for $x_i = (x,i) \in X\times\{0,1\}$ and $A\in \mathcal{B}\times\{0,1\}$ as follows.
\begin{gather*}
    \hat{P}(x_i, A) = \begin{cases}
    P^*(x, A) & x_i \in X_0\backslash U_0\\
    \frac{2P^*(x,A)-\delta \mu^*(A)}{2-\delta} & x_i \in U_0\\
    \mu^*(A) & x_i \in X_1.
    \end{cases}
\end{gather*}
An absolutely vital aspect of this construction is that our original chain is exactly the marginal of this constructed chain and thus $\tau_{A_1\cup A_0}$ for the split chain is exactly equal to $\tau_{A}$ for the original chain. Another vital aspect is that $U_1$ now forms an atom for this new Markov chain. For more information on the splitting method in the context of Markov Chains, we refer the reader to Chapter 5 of \cite{meyn2012markov}.
\par Another important aspect of the construction is the following estimate.
\begin{align*}
    x_i \in U_1 \implies& \hat{P}(x_i, U_1) = \mu^*(U_1) = \frac{\delta}{2}\\
    x_i \in U_0 \implies& \hat{P}(x_i,U_1) =  \frac{2P^*(x,U_1)-\delta \mu^*(U_1)}{2-\delta}\\
    &= \frac{\delta \left(P(x,U) - \frac{\delta}{2}\mu(U)\right) }{2-\delta}
     \geq \frac{\delta}{2-\delta} \frac{\delta}{2}.
\end{align*}
In either case, we have $\mathbb{P}_{x_i}\left(\tau_{U_1} = 1\right) \geq \frac{\delta}{2-\delta} \frac{\delta}{2}$ (note this is always less than or equal to $\frac{\delta}{2}$ as necessarily $0< \delta \leq 1$). 
\par Now we have quantitative control from $U_0\cup U_1$ to $U_1$ and so we can again apply Corollary \ref{MainCor} to show that $U_1$ is an atom with $\mathbb{E}_{x}\left[r^{\tau_{U_1}}\right] < \infty$ for all $x\in X_0\cup X_1$ and bounded on $U_1$ for some effective constant $r$ with effective bounds. 
The goal is to relate this control of geometric tails on the first occupation time for an atom to a control over the $V$-weighted distance. The key tool for this part of the proof will be the so-called ``regenerative" or ``first entrance-last exit" decomposition. 
\begin{dfn}
Given a transition kernel $P(x,dy)$ generating a Markov chain $\{\Phi^n\}_{n=1}^\infty$ on a measurable space $(X,\mathcal{B})$ and a $B\in \mathcal{B}$, the \textbf{taboo probability} for the set is the probability of reaching another set while avoiding $B$. More formally, it is defined as 
$$\prescript{}{B}{P}^N(x,A) = \mathbb{P}_{x}\left\{ \Phi^n \in A, \tau_B \geq n\right\}.  $$
This can also be expressed inductively according the the following relations.
\begin{align*}
    \prescript{}{B}{P}(x,A) &= P(x,A)\\
    \prescript{}{B}{P}^{N+1}(x,A) &= \int_{B^c}\prescript{}{B}{P}^{N+1}(x,dy)P(y,A).
\end{align*}
\end{dfn}
\begin{dfn}
Given a transition kernel $P(x,dy)$ on a measurable space $(X,\mathcal{B})$ and a $B\in \mathcal{B}$, we define the \textbf{regenerative decomposition} of the kernel as follows
$$P^n(x,A) = \prescript{}{B}{P}^n(x,C) + \sum_{j=1}^{n-1}\int_B \left[\sum_{k=1}^{j} \int_B \prescript{}{B}{P}^k(x,dv)P^{j-k}(v,dw) \right] \prescript{}{B}{P}^{n-j}(w,A).$$
\end{dfn}
Intuitively, the regenerative decomposition tracks the first entrance and last exit from a given set. Every trajectory can be divided into never meeting the set in $n$ steps or meeting in less than $n$ steps. These trajectories can again be divided by how long until they have their last visit to $B$ before leaving and never returning within the first $n$ iterates.  Notice that the last exit time in general is not a stopping time. The fact that we only consider the first $n$ steps at a time is therefore crucial.
\par If the set $B = \alpha$ is an atom and we are considering the Markov chain that we have defined on the split chain, this simplifies from a series of integrals into  a sum as follows.
\begin{gather*}
    \hat{P}^n(x,C) = \prescript{}{\alpha}{\hat{P}}^n(x,C) + \sum_{j=1}^{n-1} \sum_{k=1}^{j}  \prescript{}{\alpha}{\hat{P}}^k(x,\alpha)\hat{P}^{j-k}(\alpha,\alpha) \prescript{}{\alpha}{\hat{P}}^{n-j}(\alpha,C)
\end{gather*}
We can then consider the difference with a fixed measure (say the invariant measure). We must first verify that the invariant measure exists. As we already know that $\mathbb{E}_{\alpha}[r^\tau_\alpha]$ is bounded, we have that $\mathbb{E}_{\alpha}[\tau_\alpha]$ is also finite. Thus by Theorem 10.2.1 in \cite{meyn2012markov}, there exists an invariant measure for the split chain. A key additional fact is that $\hat{\pi}(\alpha) = \mathbb{E}_{\alpha}[\tau_\alpha]^{-1}.$
\par From a direct application of Theorem 10.4.9 in \cite{meyn2012markov}, we have that
$$\pi(g) = \pi(\alpha)\mathbb{E}_{\alpha}\left[\sum_{k=0}^{\tau_\alpha - 1}g\left(\hat{\Phi}^k\right) \right] = \pi(\alpha) \sum_{k=1}^\infty \mathbb{E}_{\alpha}\left[ g\left(\hat{\Phi}^k\right) \mathbbm{1}_{\tau_\alpha \geq k}\right] = \pi(\alpha)\sum_{k=1}^\infty \prescript{}{\alpha}{\hat{P}}^k(\alpha,g) .$$
Based on these two decompositions, we introduce three sequences of value.
\begin{align*}
    a_x(n) & := \mathbb{P}_{x}\left\{\tau_\alpha = n \right\} = \prescript{}{\alpha}{\hat{P}}^n(x,\alpha)\\
    u(n) & := \mathbb{P}_{\alpha}\left\{ \hat{\Phi}^n \in \alpha \right\} = \hat{P}^n(\alpha,\alpha)\\
    t_g(n) & := \mathbb{E}_{\alpha}\left[ g\left(\hat{\Phi}^n\right) \mathbbm{1}_{\tau_\alpha\geq n} \right] = \prescript{}{\alpha}{\hat{P}}^n(\alpha,g)
\end{align*}
Combining these terms, we have the following bound.
\begin{align*}
    \hat{P}^n(x,g) & = \prescript{}{\alpha}{\hat{P}}^n(x,g) + a_x*u*t_g(n)\\
    \hat{\pi}(g) & = \hat{\pi}(\alpha) \sum_{k=1}^\infty \prescript{}{\alpha}{\hat{P}}^k(\alpha,g)= \hat{\pi}(\alpha)* t_g(n) + \hat{\pi}(\alpha)\sum_{k=n+1}^\infty \prescript{}{\alpha}{\hat{P}}^k(\alpha,g)
\end{align*}
We have used a small abuse of notation by saying that $\hat{\pi}(\alpha)*t_g(n)$ is the convolution of $t_g(n)$ with the constant sequence of value $\hat{\pi}(\alpha)$.
We then obtain the following useful bound.
\begin{align*}
    \abs{\hat{P}^n(x,g) - \hat{\pi}(g)} \leq \prescript{}{\alpha}{P}^n(x,g) + \abs{a_x*u - \hat{\pi}(\alpha)}*t_g(n) + \hat{\pi}(\alpha)\sum_{k=n+1}^\infty \prescript{}{\alpha}{P}^k(\alpha,g).
\end{align*}
Multiplying by $r^n$ and summing, this gives the following.
\begin{align*}
    \sum_{n=1}^\infty \abs{\hat{P}^n(x,g) - \hat{\pi}(g)}r^n &\leq \sum_{n=1}^\infty \prescript{}{\alpha}{\hat{P}}^n(x,g)r^n + \sum_{n=1}^\infty \abs{a_x*u - \hat{\pi}(\alpha)}*t_g(n)r^n \\
    &+ \hat{\pi}(\alpha)\sum_{n=1}^\infty \sum_{k=n+1}^\infty \prescript{}{\alpha}{\hat{P}}^k(\alpha,g)r^n.
\end{align*}
Our goal is to now bound these terms. Following the bounds in pages 360-361 of \cite{meyn2012markov}, we achieve the following for all $\|g\|_{\hat{V}} \leq 1$. 
\begin{align*}
   \sum_{n=1}^\infty \abs{\hat{P}^n(x,g) - \hat{\pi}(g)} r^n &\leq \mathbb{E}_{x}\left[\sum_{n=1}^{\tau_\alpha} \hat{V}(\Phi^n)r^n\right] + \frac{1}{r-1}\mathbb{E}_{\alpha}\left[\sum_{n=1}^{\tau_\alpha} \hat{V}(\Phi^n)r^n\right] \\
   &+ \mathbb{E}_{\alpha}[r^{\tau_{\alpha}}]\left[\frac{r}{r-1} + \sum_{n=1}^\infty \abs{u(n) - \hat{\pi}(\alpha)}r^n\right]
\end{align*}
Our goal will now be to bound these terms. 
\par First notice that
\begin{align*}
    \mathbb{E}_{x}\left[\sum_{n=1}^{\tau_\alpha} \hat{V}(\Phi^n)r^n\right] &\leq \sup_{x\in \alpha} \hat{V}(x) \mathbb{E}_{x}[r^{\tau_\alpha}] + \mathbb{E}_{x}\left[\sum_{n=0}^{\tau_\alpha-1} \hat{V}(\Phi^n)r^n\right]\\
    \mathbb{E}_{\alpha}\left[\sum_{n=1}^{\tau_\alpha} \hat{V}(\Phi^n)r^n\right] &\leq \sup_{x\in \alpha} \hat{V}(x) \mathbb{E}_{\alpha}[r^{\tau_\alpha}] + \mathbb{E}_{\alpha}\left[\sum_{n=0}^{\tau_\alpha-1} \hat{V}(\Phi^n)r^n\right]
\end{align*}
By definition, $$\sup_{x\in \alpha} \hat{V}(x) =\sup_{x\in U_1} \hat{V}(x) = \sup_{x\in U} V(x) = M_U. $$
The other terms are then bounded by effective constants multiplied by $V(x).$ 
\par We now turn to the last term. By construction, we have that $P(\alpha,\alpha) = p(1) \geq \frac{\delta}{2}$ and that $\mathbb{E}_{\alpha}\left[r^{\tau_\alpha}\right] < \infty.$ Comparing forms directly, we also see that $\pi(1) = \hat{\pi}(\alpha).$
Thus by Theorem \ref{MainKendallTheorem}, we have that there exists an effective constant $D < \infty$ such that $\sum_{n=1}^\infty \abs{u(n) - \hat{\pi}(\alpha)}r^n \leq D$. Again, we can bound $\mathbb{E}[r^{\tau_\alpha}]$ by an effective constant. Since $V(x) \geq 1$, we can multiply all relevant terms by it without loss to obtain an effective constant $D$ with the following bound for all $g$ such that $\|g\|_{\hat{V}}\leq 1$.
\begin{gather*}
    \sum_{n=1}^\infty \abs{\hat{P}^n(x,g) - \hat{\pi}(g)} r^n  \leq D V(x)\\
    \sup_{g: \|g\|_{\hat{V}}\leq 1}\abs{\hat{P}^n(x,g) - \hat{\pi}(g)} \leq D V(x) r^{-n}.
\end{gather*}
Using the fact that our original Markov chain is the marginal of this chain, we have the final bound exactly as in page 327 of \cite{meyn2012markov}. $\pi$ here indicates the marginal of $\hat{\pi}$. This is an invariant measure of the original Markov process by Theorem 10.4.1 in \cite{meyn2012markov}. This then gives the theorem.
\begin{gather*}
    \sup_{g: \|g\|_{V}\leq 1}\abs{P^n(x,g) - \pi(g)}  \leq \sup_{g: \|g\|_{\hat{V}}\leq 1}\abs{\hat{P}^n(x,g) - \hat{\pi}(g)} \leq D V(x) r^{-n}.
\end{gather*}
\printbibliography
\end{document}